\documentclass{article}
\usepackage[a4paper, margin=1in]{geometry}
\usepackage{amsmath, amssymb, amsthm}
\usepackage{graphicx}
\usepackage{hyperref}
\usepackage{amsmath}     				
\usepackage{amssymb}   					
\usepackage{amsthm}     				
\usepackage{graphicx}    				

\usepackage{textcomp}
\usepackage[utf8]{inputenc}
\usepackage[english]{babel}    				
\usepackage[T1]{fontenc} 				
\usepackage{marvosym}    				
\usepackage[sc]{mathpazo}  	 			

\usepackage{authblk} 					
\usepackage[usenames]{xcolor}  			
\usepackage{lastpage} 					

\usepackage{enumerate}	 				

\usepackage{hyperref} 					

\usepackage{amssymb}

\newcommand{\cF}{{\mathcal F}}

\newcommand{\cT}{{\mathcal T}}

\newcommand{\norm}[1]{\left\|#1\right\|}				

\newcommand{\tm}{\times}
\newcommand{\fall}{\quad\text{for all }}
\newtheorem{theorem}{Theorem}[section]

\newtheorem{lemma}[theorem]{Lemma}

\theoremstyle{definition}
\newtheorem{definition}[theorem]{Definition}
\theoremstyle{definition}
\newtheorem{remark}[theorem]{Remark}
\theoremstyle{definition}

\numberwithin{equation}{section}
\numberwithin{table}{section}
\numberwithin{figure}{section}
\title{$C^m-$ linearization of discrete random dynamical systems}
\author{Iryna Vasylieva \\
        \small University of Klagenfurt, A--9020 Universit\"atsstraße 65--67, Klagenfurt,  Austria \\
        \small \texttt{iryna.vasylieva@aau.at}}
\date{}

\begin{document}

\maketitle

\begin{abstract}
    This paper establishes $C^m$ topological equivalence of nonautonomous semilinear difference equation with its linearization and generalizes the obtained results to discrete random dynamical systems, considering both, global and local, assumptions.
\end{abstract}
\smallskip
\section{Introduction} 

The study of linearization of differential and difference equations has a rich history, originating from the classical Hartman–Grobman Theorem \cite{grobman:59,hartman:60a,hartman:60b}. As a fundamental tool in dynamical systems, linearization provides insight into the local behavior of nonlinear systems by approximating them with simpler linear models. The Hartman–Grobman theorem ensures the existence of a local homeomorphism that maps a nonlinear system to its linearized counterpart near a fixed point, provided that the associated linearized flow satisfies the hyperbolicity condition.

Since the original result, the Hartman–Grobman theorem has been extensively generalized into various directions.

Speaking about generalizations, it is important to mention results for
nonautonomous difference equations \cite{Aul_06}, where the authors extend the results of \cite{palmer:75} to the discrete-time setting. Their work demonstrates that under the assumption of an exponential dichotomy for the linear system and a small Lipschitz constant for the nonlinear terms, the Hartman-Grobman theorem remains valid.  Additionally,  \cite{Cas_20} consider  difference equations with a nonlinear perturbation and establishes  sufficient conditions for topological
equivalence between the aforementioned systems.

The smoothness qualities of the linearization map are an important research subject in the linearization theory because they define the amount of structural information that is retained during the move from a nonlinear system to the linear one. This problem has been extensively studied in deterministic settings. For example, \cite{bar_06grobman} extends the Hartman-Grobman theorem to Banach spaces for nonuniform hyperbolic dynamics, proving the existence of Hölder continuous conjugacies and analyzing their dependence on perturbations. Similarly, \cite{dragiv19} establishes $C^1$ smooth linearization for nonautonomous difference equations with a nonuniform strong exponential dichotomy, extending the results to infinite-dimensional spaces. Additionally, \cite{Cas_21} proves that for discrete nonautonomous equations, the homeomorphisms in the linearization theorem are $C^2$ diffeomorphisms, avoiding spectral gaps and non-resonance conditions.

The Hartman–Grobman theorem has been extended to random dynamical systems, which include stochasticity or randomness in their development, with growing interest in recent years (see, e.g., \cite{Arnold}, \cite{Chu_02}). These systems are especially useful for modeling real-world problems where it is important to pay attention on uncertainties or random perturbations.
Early extensions of the theorem to random dynamical systems trace back to the work of Arnold~\cite{Arnold}, building on results by Wanner~\cite{Wanner_95}.
The work  \cite{Li_05} extends Sternberg-type linearization theorems to random dynamical systems, proving a stable and unstable manifold theorem with tempered estimates. Similarly, \cite{Lu_05} generalizes classical linearization results to analytic random dynamical systems by establishing Poincaré-type and Poincaré–Dulac-type theorems. Additionally, \cite{Zhao_20} provides a simplified proof of the Hartman–Grobman theorem for random dynamical systems, eliminating the need for preexisting stable and unstable manifolds by using an iterative approach to establish topological conjugacy.

The smoothness properties of the conjugating map in random dynamical systems are as interesting as in deterministic cases. At this point, it is worth mentioning the
paper \cite{Lu_20}, which extends the $C^1$ Hartman theorem to random dynamical systems, proving $C^{1,\beta}$ linearization result.

The aim of this paper is to establish a linearization result for discrete random dynamical systems in Banach spaces, with a particular focus on determining the smooth properties of the conjugating map.

The  paper is arranged as follows. In the remaining part of this section, we introduce the basic notations and definitions needed throughout the paper. Section 2 focuses on nonautonomous semilinear difference equations, where we establish the smooth topological equivalence results. In Section 3, we apply these results to the discrete random dynamical systems, deriving global and local smooth \( C^m \)-linearization results. Some important theorems that we used in proving our results are included in the Appendix.

\subsection{Preliminaries}~\label{Prel}

Let $X = (X,\|{\cdot}\|)$ be an arbitrary Banach space and let $L(X)$ denote the space of all bounded linear operators on X. Throughout let $\mathbb{I}\subset \mathbb{Z}$ be a discrete interval unbounded above, that is the intersection of a real interval unbounded above with the set $\mathbb{Z}.$   By $I_X$ we denote the identity operator on $X.$
Moreover, let $\{\|{\cdot}\|_{t,\omega}\}_{t\in I,\,\omega\in\Omega}$ represent a family of norms on $X$ such that
\begin{enumerate}
	\item[(i)] the mapping $(t,\omega,x)\mapsto\|{x}\|_{t,\omega}$ is strongly measurable,

	\item[(ii)] there exists a measurable mapping $\ell:\mathbb{I}\times\Omega\to(0,\infty)$ such that $\ell(\cdot,\omega)$ is locally bounded  and
	\begin{equation}
		\frac{1}{\ell(t,\omega)}\|{x}\|\leq\|{x}\|_{t,\omega}\leq\ell(t,\omega)\|{x}\|
		\text{ for all } t\in \mathbb{I},\,\omega\in\Omega,\,x\in X.
		\label{normest}
	\end{equation}
\end{enumerate}
For a linear operator $T\in L(X)$ we furthermore define
$$
	\|{T}\|_{s,t,\omega}:=\sup_{\|{x}\|_{s,\omega}=1}\|{Tx}\|_{t,\omega}\text{ for all } s,t\in \mathbb{I},\,\omega\in\Omega
$$
and obtain $\|{Tx}\|_{t,\omega}\leq\|{T}\|_{s,t,\omega}\|{x}\|_{s,\omega}$ for all $t,s\in \mathbb{I}$ and $x\in X$.

Given this, a function $\mu:\mathbb{I}\to X$ is denoted as \emph{bounded} (for $\omega\in\Omega$), if
$$
	\|{\mu}\|_\omega:=\sup_{t\in \mathbb{I}}\|{\mu(t)}\|_{t,\omega}<\infty.
$$

Let  $(\Omega, \mathcal{F})$ be an arbitrary measurable space, and
let $A:\mathbb{I}\times\Omega\to GL(X)$  be a strongly measurable map. Consider the linear difference equation
\begin{equation}\label{lin1}
x_{t+1} = A_{\omega}(t)x_t.
\end{equation}

The general solution of \eqref{lin1} is given by the evolution operator
$\Phi_{\omega}(t,s)$ in the following way
\begin{equation}
    \Phi_{\omega}(t,s) = \begin{cases}
   I_{X} &t = s,\\
   A_{\omega}(t-1)\circ...\circ A_{\omega}(s) &s < t,\\
   A_{\omega}(s)^{-1}\circ...\circ A_{\omega}(t-1)^{-1} &s > t,
 \end{cases}
\end{equation}
for arbitrary $s,t\in\mathbb{I}$ and $\omega\in\Omega.$
And it has the properties
$$\Phi_{\omega}(t+1,s) = A_{\omega}(t)\circ\Phi_{\omega}(t,s),$$
$$\Phi_{\omega}(t,s) = \Phi_{\omega}(t,p)\circ\Phi_{\omega}(p,s),$$
$$\Phi_{\omega}(t,s)^{-1} = \Phi_{\omega}(s,t),$$
for all $s,p\in\mathbb{I},$ $t\in\mathbb{I}$ and $\omega\in\Omega.$

\begin{definition}
 Let $(\Omega, \mathcal{F})$ be an arbitrary measurable space, let $f:\mathbb{I}\times X\times \Omega\to X$ be an arbitrary mapping . Then the equation
 \begin{equation}\label{par_eq}
   x_{t+1} = f(t,x_t,\omega)
 \end{equation}
 is called a \emph{nonautonomous parameter dependent difference equation} and the sequence $\phi = (\phi(t))_{t\in\mathbb{I}}$ is called a \emph{solution} of \eqref{par_eq} (to the parameter $\omega$), if the  identity
 \begin{equation}\label{sol_par}
  \phi(t+1) = f(t,\phi(t),\omega)
 \end{equation}
 holds for all $ t\in \mathbb{I}.$
 If in addition the mapping $\phi$ satisfies the condition
 \begin{equation}\label{ivp}
    x_{\tau} = \xi,
 \end{equation}
 i.e. $\phi(\tau) = \xi,$ then we say that $\phi$ solves initial value problem \eqref{par_eq}, \eqref{ivp} for $\tau\in\mathbb{I},$ $\xi\in X$.

 \end{definition}

Obviously, for given $\tau\in\mathbb{I}, \xi\in X,\omega\in\Omega$ the solution of the initial value problem  \eqref{par_eq}, \eqref{ivp} exists on the discrete interval $\{t\in\mathbb{I}:t\geq \tau\}$ and is uniquely determined there. This solution can be given recursively as follows
\begin{equation}
    \phi(t,\tau,\omega,\xi) = \begin{cases}
      \xi, \text{ for } t = \tau\\
      f(t-1,\phi(t-1,\tau,\omega,\xi),\omega),\text{ for } t > \tau
    \end{cases}.
\end{equation}

If additionally the mapping $f(t,\cdot,\omega)$ is bijective for any $t\in\mathbb{I}$ and $\omega\in\Omega$, then there exists a uniquely determined solution of the  initial value problem  \eqref{par_eq}, \eqref{ivp} on the whole $\mathbb{I},$ namely
\begin{equation}\label{def_gen}
    \phi(t,\tau,\omega,\xi) = \begin{cases}
      f^{-1}(t,\phi(t+1,\tau,\omega,\xi),\omega),\text{ for } t < \tau\\\xi, \text{ for } t = \tau\\
      f(t-1,\phi(t-1,\tau,\omega,\xi),\omega),\text{ for } t > \tau
    \end{cases}.
\end{equation}

In both cases, the mapping $\phi$ is called the general solution of the difference equation \eqref{par_eq}.

It is easy to verify that the general solution satisfy the identity $$\phi(t,\tau,\omega,\xi) = \phi(t,s,\omega,\phi(s,\tau,\omega,\xi))$$
for any $t,\tau,s\in \mathbb{I},$ $\omega\in\Omega$ and $\xi\in X,$ where in non-invertible case $\tau\leq s\leq t$ must be assumed.

\section{Semilinear difference equations}

In this section, we study semilinear difference equations with parameters, focusing on their topological equivalence to linear systems. Our goal is to establish
$C^m-$ linearization results under appropriate conditions.

 Let $\mathbb{I}\subset \mathbb{Z}$ with $\mathbb{I}\neq \emptyset,$ let $(\Omega, \mathcal{F})$ be an arbitrary measurable space, $F:\mathbb{I}\times\Omega\times X\to X$ and  $A:\mathbb{I}\times\Omega\to GL(X).$

Consider the semilinear difference equation
\begin{equation}\tag{N}
x_{t+1} = A_{\omega}(t)x_t+F_{\omega}(t,x_t).\label{deq}
\end{equation}

Together with the system \eqref{deq} consider a linear system
\begin{equation}\tag{$L$}
   x_{t+1} = A_{\omega}(t)x_t.
    \label{lin}
\end{equation}

Assume that the following assumptions hold for all $\omega\in\Omega$
\begin{enumerate}
	\item[$(H_1)$] the linear system \eqref{lin} has the property of bounded growth, i.e. there exist reals $K(\omega)\geq 1$ and $\alpha\in (0,1)$ such that	
	\begin{equation}
		\|{\Phi_\omega(t,s)}\|_{s,t,\omega}\leq K(\omega){\alpha}^{t-s}\text{ for all } t,s\in \mathbb{I},\,s\leq t.
		\label{hyp11}
	\end{equation}

	\item[$(H_2^0)$] there exist reals $L(\omega), M(\omega)\geq 0$ such that
	\begin{align}
		F_\omega(t,0)&=0,
		\label{hyp20}\\
		\|F_\omega(t,x)\|_{t+1,\omega}&\leq M(\omega), \label{hyp21}\\
  \|F_\omega(t,x)-F_\omega(t,\bar x)\|_{t+1,\omega}&\leq L(\omega)\|x-\bar x\|_{t,\omega}
	\label{hyp22}
	\end{align}
 for  all $t\in \mathbb{I}$ and all $x,\bar x\in X.$
\end{enumerate}
\begin{lemma}
    Under the above hypothesis, with $\displaystyle L(\omega)\|A_{\omega}(t)^{-1}\|_{t,t+1,\omega}<{1},$ the operator $A_{\omega}(t)+F_{\omega}(t,\cdot)$ is invertible for any $t\in\mathbb{I}$.
\end{lemma}

\begin{proof}
Let us define the new operator $ B_{\omega}(t)$ as follows
$$
B_{\omega}(t) := A_{\omega}(t) + F_{\omega}(t, \cdot).
$$
We need to show that for any $y\in X,$ there exists a unique $x\in X$ such that:
$B_{\omega}(t) x = y,$
or equivalently,
$A_{\omega}(t) x + F_{\omega}(t, x) = y.$
Since $ A_{\omega}(t)\in GL(X),$  we obtain
$x + A_{\omega}(t)^{-1} F_{\omega}(t, x) = A_{\omega}(t)^{-1} y.$

Define the new operator
$G(x) := A_{\omega}(t)^{-1} F_{\omega}(t, x),$
so that the equation becomes
$$x = A_{\omega}(t)^{-1} y - G(x).$$

We now show that $G$ is a contraction on X, Indeed,
\begin{eqnarray*}
    \|G(x_1) - G(x_2)\|_{t,\omega}& = &\|A_{\omega}(t)^{-1}  F_{\omega}(t, x_1) - A_{\omega}(t)^{-1}  F_{\omega}(t, x_2)\|_{t,\omega}\\
    &\leq& L(\omega)\|A_{\omega}(t)^{-1}\|_{t,t+1,\omega}\| x_1 - x_2\|.
\end{eqnarray*}
Lipschitz constant $L(\omega)$ can be choosen small enough such that
$\| A_{\omega}(t)^{-1} \|_{t,t+1,\omega} L(\omega) < 1.$
This implies that the operator $G$  is a contraction on $X$, therefore  the equation
above
has a unique fixed point by the Banach fixed-point theorem. Consequently,   $ A_{\omega}(t) + F_{\omega}(t, \cdot)$ is invertible.

\end{proof}

Before presenting the linearization results, we will first introduce some technical lemmas.

\begin{lemma}
	Let $\tau\in \mathbb{I},$ $l_{\omega}$ denotes the space of all sequence bounded with respect to $\|\cdot\|_{t,\omega}$  and $(H_1)$ holds. Then
	\begin{align}
		\mathcal L_{\tau,\omega}:X&\to l_\omega,&
		\mathcal L_{\tau,\omega}\xi&:=\Phi_\omega(\cdot,\tau)\xi
	\end{align}
	defines a sequence of bounded linear operator for all  $\omega\in\Omega$ and $\xi\in X.$
\end{lemma}\label{leml}

\begin{proof}
Let us fix $t\in\mathbb{I},$ we obtain
\begin{align}\label{norm_cL}
		\|{[\mathcal L_{\tau,\omega}\xi](t)}\|_{\tau,\omega}
		&=
		\|\Phi_\omega(t,\tau)\xi\|_{\tau,\omega}
		\stackrel{\eqref{hyp11}}{\leq}
		K(\omega)\alpha^{(t-\tau)}\|{\xi}\|_{\tau,\omega}
		=
		K(\omega)\alpha^{(t-\tau_0)}\alpha^{(\tau_0-\tau)}\|{\xi}\|_{\tau,\omega}\nonumber\\
		&\leq
		K(\omega)\alpha^{(\tau_0-\tau)}\|{\xi}\|_{\tau,\omega}
		\stackrel{\eqref{normest}}{\leq}
		K(\omega)\alpha^{(\tau_0-\tau)}\ell(\tau,\omega)\|{\xi}\|\text{ for all }\tau\leq\tau_0\leq t
	\end{align}
	and passing to the supremum over $t\in \mathbb{I}$ yields $\|{\mathcal{L}_{\tau,\omega}\xi}\|_\omega\leq K(\omega)\alpha^{(\tau_0-\tau)}\ell(\tau,\omega)\|{\xi}\|$, which is the claim.
 \end{proof}
 \begin{lemma}
	\label{lemf}
	If $(H_1)$ and $(H_2^0)$ hold on $\mathbb{I}$, then the nonlinear operator
	\begin{align}
		\mathcal{F}_\omega:l_{\omega}(\mathbb{I},X)&\to l_{\omega}(\mathbb{I},X) ,&
		\mathcal{F}_\omega(\phi)(t)&:=\sum\limits_{s=\tau_0}^{t-1}\Phi_\omega(t,s+1)F_\omega(s,\phi(s))
	\end{align}
	is well-defined and for all $\phi,\bar\phi\in X,\,\omega\in\Omega$ satisfies
	\begin{equation}\label{eq2.8}
		\|{\mathcal{F}_\omega(\phi)-\mathcal{F}_\omega(\bar\phi)}\|_\omega
		\leq
		\frac{K(\omega)L(\omega)}{1-\alpha}\|{\phi-\bar\phi}\|_\omega,
	\end{equation}
 \begin{equation}\label{eq2.9}
		\|{\mathcal{F}_\omega(\phi)}\|_\omega
		\leq
		 \frac{K(\omega)M(\omega)}{1-\alpha}.
	\end{equation}
\end{lemma}
\begin{proof}
In order to show boundness of the operator $\mathcal{F}_{\omega}$ let us fix $t\in\mathbb{I}$ and consider 	
\begin{align*}
		\|{\mathcal{F}_\omega(\phi)(t)}\|_{t+1,\omega}
		&{\leq} \sum\limits_{s=\tau_0}^{t-1}\|{\Phi_\omega(t,s+1)}\|_{s,t,\omega}\|{F_{\omega}(s,\phi(s))}\|_{s+1,\omega}\\
        &\stackrel{\eqref{hyp11}}{\leq} K(\omega) \sum\limits_{s=\tau_0}^{t-1}\alpha^{t-s}\|{F_{\omega}(s,\phi(s))}\|_{s+1,\omega}\\&
         \stackrel{\eqref{hyp21}}{\leq} K(\omega) M(\omega) \sum\limits_{s=\tau_0}^{t-1}\alpha^{t-s} \\&\leq K(\omega)M(\omega) \frac{\alpha^{t-\tau_0+1}-\alpha}{\alpha-1}\text{ for all }\tau_0\leq t
	\end{align*}
and passing to the supremum over $t\in \mathbb{I}$ yields \eqref{eq2.9}.

Moreover,
\begin{align*}
		\|{\mathcal{F}_\omega(\phi)(t)-\mathcal{F}_\omega(\bar\phi)(t)}\|_{t+1,\omega}
		&{\leq} \sum\limits_{s=\tau_0}^{t-1}\|{\Phi_\omega(t,\tau)}\|_{s,t,\omega}\|{F_{\omega}(s,\phi(s))- F_{\omega}(s,\bar\phi(s))}\|_{s+1,\omega} \\
        &{\stackrel{\eqref{hyp11}}{\leq}} K(\omega) \sum\limits_{s=\tau_0}^{t-1}\alpha^{t-s}\|{F_{\omega}(s,\phi(s))- F_{\omega}(s,\bar\phi(s))}\|_{s+1,\omega}
        \\ &{\stackrel{\eqref{hyp22}}{\leq}}K(\omega) L(\omega) \sum\limits_{s=\tau_0}^{t-1}\alpha^{t-s}\|{\phi(s)-\bar\phi(s)}\|_{s,\omega}
        \\ &{\leq} K(\omega)L(\omega)\frac{\alpha^{t-\tau_0+1}-\alpha}{\alpha-1}\|{\phi-\bar\phi}\|_{\omega}\text{ for all }\tau_0\leq t.
	\end{align*}
Passing to the supremum over $t\in \mathbb{I}$ yields \eqref{eq2.8} for all \mbox{$\phi,\bar\phi\in X$}.
\end{proof}

\begin{theorem}\label{thm_lin}

	If $(H_1)$ and $(H_2^0)$ hold on $\mathbb{I}$ with
	$$
		\displaystyle
K(\omega)L(\omega)<1-\alpha,
	$$
	then there exists a   map $H_{\omega}:\mathbb{I}\times X\to X$ such that for every $\omega\in\Omega,$ $s,t\in \mathbb{I}$
\begin{itemize}
  \item[(i)] $H_\omega(t,\cdot)$ is a homeomorphism on $X$ with $H_{\omega}(t,0) = 0;$
  \item[(ii)] the linear equation \eqref{lin} and the semilinear equation \eqref{deq} are conjugated in the sense that
  \begin{equation*}
    H_\omega(t,\Phi_{\omega}(t,s)\cdot) =\varphi_{\omega}(t,s, H_{\omega}(s,\cdot));
  \end{equation*}
   \item[(iii)] the operator $H_{\omega}(t,\cdot)$ and its inverse are near identity in the sense that for all $\xi\in X,$ $\eta\in X$:
   $$\|H_{\omega}(t,\xi)-\xi\|_{t+1,\omega}\leq \frac{K(\omega)M(\omega)}{1-\alpha}$$ and $$\|{H_{\omega}(t,\cdot)^{-1}(\eta)-\eta}\|_{t+1,\omega}\leq \frac{K(\omega)M(\omega)}{1-\alpha};$$
   \item[(iv)]  the operator $H_{\omega}(t,\cdot):X\to X$ and its inverse satisfy a global Lipschitz condition.

\end{itemize}
\end{theorem}
\begin{proof}

\emph{Step 1.}

	We fix $\tau\in \mathbb{I}, \omega\in \Omega.$ Let us define  $\mathcal{T}_{\omega}:{l}_\omega(\mathbb{I},X)\times X\to {l}_\omega(\mathbb{I},X)$,
	\begin{equation}\label{oper_tau}
		\mathcal{T}_{\omega}(\phi;\xi,\tau):=\mathcal{F}_\omega(\phi+\mathcal{L}_{\tau,\omega}\xi),
	\end{equation}
where $\mathcal{F}_\omega:l_\omega(\mathbb{I},X)\to l_\omega(\mathbb{I},X)$ is introduced in Lemma \ref{lemf} and $\mathcal{L}_{\tau,\omega}:X\to l_\omega(\mathbb{I},X)$ is introduced in Lemma \ref{leml}. We have

\begin{align*}
  \|{\mathcal{T}_{\omega}(\phi;\xi,\tau)-\mathcal{T}_{\omega}(\bar{\phi};\xi,\tau)}\|_{\omega}&\overset{\eqref{oper_tau}}{=}\|{\mathcal{F}_\omega(\phi+\mathcal{L}_{\tau,\omega}\xi)-\mathcal{F}_\omega(\bar{\phi}+\mathcal{L}_{\tau,\omega}\xi)}\|_{\omega}\\
		&{\stackrel{\eqref{eq2.8}}{\leq}}\frac{K(\omega)L(\omega)}{1-\alpha}\|{\phi-\bar{\phi}}\|_{\omega}
\end{align*}
for all \mbox{$\phi,\bar\phi\in l_\omega(\mathbb{I},X).$}
 Hence, $\mathcal{T}_{\omega}$ is a contraction uniformly in $\xi\in X$ and in $\tau\in \mathbb{I}$ and by uniform contraction principle there is a unique fixed point $\varphi^*_{\tau,\omega}(\xi)\in l_{\omega}(\mathbb{I},X)$ of $\mathcal{T}_{\omega}(\cdot;\xi,\tau)$ for all $\xi\in X,$ $\tau\in \mathbb{I}$.

 Moreover,
 \begin{align}
  \|{\mathcal{T}_{\omega}(\phi;\xi,\tau)-\mathcal{T}_{\omega}(\phi;\bar{\xi},\tau)}\|_{\tau+1,\omega}&\overset{\eqref{oper_tau}}{=}\|{\mathcal{F}_\omega(\phi+\mathcal{L}_{\tau,\omega}\xi)-\mathcal{F}_\omega(\phi+\mathcal{L}_{\tau,\omega}\bar{\xi})}\|_{\tau+1,\omega}\nonumber\\
		&{\stackrel{\eqref{eq2.8}}{\leq}}\frac{K(\omega)L(\omega)}{1-\alpha}\|{\mathcal{L}_{\tau,\omega}\xi-\mathcal{L}_{\tau,\omega}\bar{\xi}}\|_{\tau,\omega}\nonumber\\&{\stackrel{\eqref{norm_cL}}{\leq}}\frac{K^2(\omega)L(\omega)}{1-\alpha}\alpha^{(\tau_0-\tau)}\ell(\tau,\omega)\|{\xi-\bar{\xi}}\|.
\end{align}

This leads us to conclude that $\mathcal{T}_{\omega}$ is Lipschitz with Lipschitz constant\\ $\displaystyle \frac{K^2(\omega)L(\omega)}{\alpha-1}\ell(\tau,\omega).$ Therefore,  the mapping $\xi\mapsto\varphi^*_{\tau,\omega}(\xi)$ is also  Lipschitz with Lipschitz constant $\displaystyle\frac{K^2(\omega)L(\omega)\ell(\tau,\omega)}{1-K(\omega)L(\omega)}.$

Furthermore, for all $\xi\in X$ and $\tau, r\in \mathbb{I},$ $\varphi^*_{\tau,\omega}(\xi)$ satisfies   the identity
\begin{equation}\label{varphi}
    \varphi^*_{\tau,\omega}(\xi)(t) = \varphi^*_{r,\omega}(\Phi_{\omega}(r,\tau)\xi)(t).
\end{equation}

Indeed,
$$\varphi^*_{\tau,\omega}(\xi)(t) = \sum\limits_{s=\tau_0}^{t-1}\Phi_{\omega}(t,s+1)F_{\omega}(s,\varphi_{\tau,\omega}^*(\xi)(s)+\Phi_{\omega}(s,\tau)\xi)$$
is a  unique solution to the initial value problem $ x_{t+1} = A_{\omega}(t)x_t + F_{\omega}(t, x_t + \Phi_{\omega}(s,\tau) \xi),
x_{\tau} = 0.
$
On the other hand,
\begin{eqnarray*}
    \varphi^*_{r,\omega}(\Phi_{\omega}(r,\tau)\xi)(t) &=& \sum\limits_{s=\tau_0}^{t-1}\Phi_{\omega}(t,s+1)F_{\omega}(s,\varphi^*_{r,\omega}(\Phi_{\omega}(r,\tau)\xi)(s)+\Phi_{\omega}(s,r)\Phi_{\omega}(r,\tau)\xi)\
    \\&=& \sum\limits_{s=r}^{t-1}\Phi_{\omega}(t,s+1)F_{\omega}(s,\varphi^*_{r,\omega}(\Phi_{\omega}(r,\tau)\xi)(s)+\Phi_{\omega}(s,\tau)\xi)
\end{eqnarray*}
is a  unique solution to the same initial value problem. As consequence of uniqueness of solution we obtain the claim.

\emph{Step 2.} For any fixed $t\in\mathbb{I},$ let us construct the maps  $H_{\omega}(t,\cdot): X\to l_{\omega}(I,X)$ and  $G_{\omega}(t,\cdot): X\to l_{\omega}(I,X)$ as follows
 \begin{align}\label{H}H_\omega{(t,\xi)}&{:=}\xi+\sum\limits_{s=\tau_0}^{t-1} \Phi_{\omega}(t,s+1)F_{\omega}(s,\varphi^*_{t,\omega}(\xi)(s)+\Phi_{\omega}(s,t)\xi)\nonumber\\
 &{=}\xi+ \varphi^*_{t,\omega}(\xi)(t).
 \end{align}
 \begin{align*}G_\omega{(t,\eta)}&{:=}\eta - \sum\limits_{s=\tau_0}^{t-1} \Phi_{\omega}(t,s+1)F_{\omega}(s,\varphi_{\omega}(s,t,\eta)).
 \end{align*}

 We can show the boundedness of $H_\omega{(t,\xi)} - \xi$ and $G_\omega{(t,\eta)} -\eta$ as follows
 \begin{align*}
     \|{H_\omega{(t,\xi)} - \xi}\|_{t+1,\omega}&{\leq}\sum\limits_{s=\tau_0}^{t-1}\|{\Phi_{\omega}(t,s+1)}\|_{s,t,\omega}\|{F_{\omega}(s,\varphi^*_{t,\omega}(\xi)(s)+\Phi_{\omega}(s,t)\xi)}\|_{s+1,\omega} \\
     &{\stackrel{\eqref{hyp11}}{\leq}}K(\omega)\sum\limits_{s=\tau_0}^{t-1} \alpha^{(t-s-1)}\|{F_{\omega}(s,\varphi^*_{t,\omega}(\xi)(s)+\Phi_{\omega}(s,t)\xi)}\|_{s+1,\omega} \\
     &{\stackrel{\eqref{hyp21}}{\leq}} K(\omega)M(\omega)\sum\limits_{s=\tau_0}^{t-1} \alpha^{(t-s-1)} \leq \frac{K(\omega)M(\omega)}{1-\alpha}\text{ for all }\xi\in X,
 \end{align*}
 and respectively
 \begin{align*}
    \|{G_\omega{(t,\eta)} - \eta}\|_{t+1,\omega}&{\leq}\sum\limits_{s=\tau_0}^{t-1}\|{\Phi_{\omega}(t,s+1)}\|_{s,t,\omega}\|{F_{\omega}(s,\varphi_{\omega}(s,t,\eta)}\|_{s+1,\omega} \\
     &{\leq}K(\omega)M(\omega)\sum\limits_{s=\tau_0}^{t-1} \alpha^{(t-s-1)}\leq \frac{K(\omega)M(\omega)}{1-\alpha}\text{ for all }\eta\in X.
 \end{align*}
The next our goal is to show that  $H_{\omega}$ maps solutions of \eqref{lin} to solutions of \eqref{deq} and $G_{\omega}$ maps solutions of \eqref{deq}
to solutions of \eqref{lin}. Using \eqref{varphi} we obtain
\begin{eqnarray*}
        H_{\omega}(t,\Phi_{\omega}(t,\tau)\xi) &=&\Phi_{\omega}(t,\tau)\xi + \sum\limits_{s=\tau_0}^{t-1} \Phi_{\omega}(t,s+1)F_{\omega}(s,\varphi_{t,\omega}^{*}(\Phi_{\omega}(t,\tau)\xi)(s)+ \Phi_{\omega}(s,\tau)\xi)\\
        &=&\Phi_{\omega}(t,\tau)\xi+\varphi_{\tau,\omega}^{*}(\xi)(t).\end{eqnarray*}
On the other side
$$H_{\omega}(t,\Phi_{\omega}(t,\tau)\xi) = \Phi_{\omega}(t,\tau)\xi+\sum\limits_{s=\tau_0}^{t-1} \Phi_{\omega}(t,s+1)F_{\omega}(s,H_{\omega}(t,\Phi_{\omega}(s,\tau)\xi)).$$
Combining both equalities above and the fact that $\Phi_{\omega}(t+1,s) = A_{\omega}(t)\circ\Phi_{\omega}(t,s),$ we obtain
        \begin{eqnarray*}
        H_{\omega}(t+1,\Phi_{\omega}(t+1,s)\xi) &=& \Phi_{\omega}(t+1,s)\xi+\varphi_{s,\omega}^{*}(\xi)(t+1)\\
        &=& A_{\omega}(t)\Phi_{\omega}(t,s)\xi+A_{\omega}(t)\varphi_{s,\omega}^{*}(\xi)(t)\\&+&F_{\omega}(t,\varphi_{s,\omega}^{*}(\xi)(t)+\Phi_{\omega}(t,s)\xi)
\\&=&A_\omega(t)H_\omega(t,\Phi_{\omega}(t,s)\xi)+F_\omega(t, H_\omega(t,\Phi_{\omega}(t,s)\xi).
    \end{eqnarray*}
        This means $t\mapsto H_{\omega}(t,\Phi_{\omega}(t,\tau)\xi)$ solves the initial value problem $x_{t+1} = A_{\omega}(t)x_t+F_{\omega}(t,x_t),$
        $x_{\tau} = H_\omega (\tau,\xi).$

As consequence of uniqueness of solution we obtain $$H_{\omega}(t,\Phi_{\omega}(t,\tau)\xi) =\varphi_{\omega}(t,\tau,H_{\omega}(\tau,\xi))$$ for all $\xi\in X.$

 In the same way it can be shown that
\begin{equation*}
  G_{\omega}(t,\varphi_{\omega}(t,s,\eta)) = \Phi_{\omega}(t,s)G_{\omega}(s,\eta) \text{ for all }\eta\in X.
\end{equation*}

\emph{Step 3.} In order to complete the proof and to show that $H_{\omega}(t,\cdot)$ is a homeomorphism we start with the claim that $H_{\omega}(t,\cdot):X\to X$ is bijective. Indeed, for any fixed $t\in \mathbb{I}$ and any $\eta\in X$ we have

\begin{align*}
     H_\omega(t,G_{\omega}(t,\varphi_{\omega}(s,t,\eta))) &\overset{\eqref{H}}{=}G_{\omega}(t,\varphi_{\omega}(s,t,\eta))
     \\
     &{+}\sum\limits_{s=\tau_0}^{t-1}\Phi_{\omega}(t,s+1)F_{\omega}(s,\varphi^*_{t,\omega}(G_{\omega}(s,\varphi_{\omega}(s,t,\eta)))\\
     &{+}\Phi_{\omega}(s,t)G_{\omega}(s,\varphi_{\omega}(s,t,\eta))))\\
     &{=}\varphi_{\omega}(s,t,\eta)-\sum\limits_{s=\tau_0}^{t-1} \Phi_{\omega}(t,s+1)F_{\omega}(s,\varphi_{\omega}(s,t,\eta))\\
     &{+}\sum\limits_{s=\tau_0}^{t-1}\Phi_{\omega}(t,s+1)F_{\omega}(s,\varphi^*_{t,\omega}(G_{\omega}(s,\varphi_{\omega}(s,t,\eta)))\\
     &{+}
     \Phi_{\omega}(s,t)G_{\omega}(s,\varphi_{\omega}(s,t,\eta))))
 \end{align*}
 If we abbreviate by $\nu(t) = \|{H_\omega(t,G_{\omega}(t,\varphi_{\omega}(s,t,\eta)))  - \varphi_{\omega}(s,t,\eta)}\|_{t+1,\omega},$ then for all  $t\in \mathbb{I}$ we have
\begin{align*}
    \nu(t) &\leq \sum\limits_{s=\tau_0}^{t-1} \|\Phi_{\omega}(t,s+1)\|_{s,t,\omega}
    \|F_{\omega}(s,\varphi^*_{t,\omega}(G_{\omega}(s,\varphi_{\omega}(s,t,\eta))) \\
    &+ \Phi_{\omega}(s,t)G_{\omega}(s,\varphi_{\omega}(s,t,\eta))
    - F_{\omega}(s,\varphi_{\omega}(s,t,\eta))\|_{s+1,\omega} \\
    &\leq K(\omega)L(\omega)\sum\limits_{s=\tau_0}^{t-1}\alpha^{(t-s-1)}
    \|\varphi^*_{t,\omega}(G_{\omega}(s,\varphi_{\omega}(s,t,\eta))) \\
    & + \Phi_{\omega}(s,t)G_{\omega}(s,\varphi_{\omega}(s,t,\eta))
    - \varphi_{\omega}(s,t,\eta)\|_{s,\omega} \\
    &\leq K(\omega)L(\omega)\sum\limits_{s=\tau_0}^{t-1}\alpha^{(t-s-1)}
    \|H_{\omega}(s,\Phi_{\omega}(s,t)G_{\omega}(s,\varphi_{\omega}(s,t,\eta))) - \varphi_{\omega}(s,t,\eta)\|_{s,\omega} \\
    &\leq K(\omega)L(\omega)\sum\limits_{s=\tau_0}^{t-1}\alpha^{(t-s-1)}\nu(s)
    \leq \frac{K(\omega)L(\omega)}{1-\alpha} \sup_{s \in \mathbb{I}'} \nu(s).
\end{align*}

From the above estimate and the assumption  $\displaystyle
K(\omega)L(\omega)<1-\alpha$ it follows that $\nu(t) = 0$ for any $t\in \mathbb{I}$ and with $t=\tau,$ $H_{\omega}(\tau,G_{\omega}(\tau,\eta)) = \eta.$

Now we prove that for any fixed $t\in \mathbb{I}$ and all $\xi\in X$ holds $G_{\omega
}(t,H_{\omega}(t,\xi)) = \xi.$ We have

\begin{align*}
    G_{\omega}(t,H_{\omega}(t,\Phi_{\omega}(t,\tau)\xi))
    &= H_{\omega}(t,\Phi_{\omega}(t,\tau)\xi)\\
    &- \sum_{s=\tau_0}^{t-1}\Phi_{\omega}(t,s+1)
    F_{\omega}(s,\varphi_{\omega}(s,t,H_{\omega}(t,\Phi_{\omega}(t,\tau)\xi)))\\
    &= \Phi_{\omega}(t,\tau)\xi
    + \sum_{s=\tau_0}^{t-1}\Big(
        F_{\omega}(s,H_{\omega}(s,\Phi_{\omega}(s,\tau)\xi)) \\
    & - F_{\omega}(s,\varphi_{\omega}(s,\tau,H_{\omega}(\tau,\xi)))
    \Big)\\
    &= \Phi_{\omega}(t,\tau)\xi.
\end{align*}

Now, if we set $t=\tau,$ then $G_{\omega}(\tau, H_{\omega}(\tau,\xi)) = \xi$ results.

As a consequence, for any $t\in \mathbb{I}$, $H_\omega(t,\cdot): X \to X$ is a bijection, and $G_{\omega}(t,\cdot): X\to X$ is its inverse.

\emph{Step 4.}
In this step  we  show that the conjugacy and its inverse satisfies Lipschitz condition with Lipschitz constants $L_{H}$ and $L_{G}$ respectively, where
$$L_{G} = 1+\frac{K^2(\omega)L(\omega)}{1-\alpha^2}e^{\frac{K(\omega)L(\omega)}{1-\alpha}}  ,\quad\quad L_{H} = 1+\frac{K^2(\omega)L(\omega)}{1-\alpha}+\frac{K^3(\omega)L^2(\omega)\alpha^{\tau_0-\tau}\ell(\tau,\omega)}{(1-\alpha)(1-K(\omega)L(\omega))}.$$

In order to show the claim for the operator $G_\omega$ we start with additional estimates.

The general solution $\varphi_{\omega}$ of \eqref{deq}  also solves the inhomogeneous equation $$x_{t+1} = A_{\omega}(t)x_t+F_{\omega}(t,\varphi_{\omega}(\tau,t,\eta)).$$
Hence, by variation of constants formula (Theorem~\ref{var_const}), we have
$$\varphi_{\omega}(s,t,\eta) = \Phi_{\omega}(t,s)\eta+\sum\limits_{\tau=s}^{t-1}\Phi_{\omega}(t,\tau+1)F_{\omega}(\tau,\varphi_{\omega}(\tau,t,\eta))$$
and we obtain
 \begin{align*}
\|{\varphi_{\omega}(s,t,\eta)-\varphi_{\omega}(s,t,\bar{\eta})}\|_{t+1,\omega}&{\leq}\|{\Phi_{\omega}(t,s)\|_{s,t,\omega}\|\eta-\bar{\eta}}\|\nonumber\\  &+\sum\limits_{\tau=s}^{t-1}\|{\Phi_{\omega}(t,\tau)}\|_{\tau,t,\omega}\|{F_{\omega}(\tau,\varphi_{\omega}(\tau,t,\eta))-F_{\omega}(\tau,\varphi_{\omega}(\tau,t,\bar{\eta}))}\|_{\tau+1,\omega} \nonumber\\&\leq
K(\omega)\alpha^{(t-s)}\|{\eta-\bar{\eta}}\|\\&+K(\omega)L(\omega) \sum\limits_{\tau=s}^{t-1}\alpha^{(t-\tau-1)}
     \|{\varphi_{\omega}(\tau,t,\eta)-\varphi_{\omega}(\tau,t,\bar{\eta})}\|_{\tau,\omega}
     \\&\leq
K(\omega)\alpha^{(t-s)}\|{\eta-\bar{\eta}}\|\\&+K(\omega)L(\omega) \sum\limits_{\tau=s}^{t-1}\alpha^{(t-\tau-1)}
     \|{\varphi_{\omega}(\tau,t,\eta)-\varphi_{\omega}(\tau,t,\bar{\eta})}\|_{\tau,\omega}.
 \end{align*}
Multiplying both sides by
 $\alpha^{(s-t)}$ and using discrete Gronwall's
inequality, gives us
\begin{align*}
\alpha^{(s-t)}\|{\varphi_{\omega}(s,t,\eta)-\varphi_{\omega}(s,t,\bar{\eta})}\|_{t+1,\omega}&{\leq} K(\omega) \left(1+\frac{K(\omega)L(\omega)}{\alpha}\right)^{t-s}\|{\eta-\bar{\eta}}\|.
\end{align*}
From the estimate above we obtain
\begin{align}\label{Gron}
\|{\varphi_{\omega}(s,t,\eta)-\varphi_{\omega}(s,t,\bar{\eta})}\|_{t+1,\omega}&{\leq} K(\omega) \left(\alpha+{K(\omega)L(\omega)}\right)^{t-s}\|{\eta-\bar{\eta}}\|.
\end{align}

    Using \eqref{Gron}, we have that
       \begin{align}\label{Lip_G}
            \|{G_\omega(t,\eta)-G_\omega(t,\bar{\eta})}\|_{t+1,\omega}&
            \leq\|{\eta-\bar{\eta}}\| + K(\omega)L(\omega)\sum\limits_{s=\tau_0}^{t-1} \alpha^{(t- s-1)}\|{\varphi_{\omega}( s, t,\bar{\eta})-\varphi_{\omega}(s, t, \eta)}\|_{s,\omega}\nonumber\\
            &\leq \|{\eta-\bar{\eta}}\|\left(1+K^2(\omega)L(\omega)\sum\limits_{s=\tau_0}^{t-1} \alpha^{t-s-1}(\alpha+KL)^{t-s}\right)\nonumber\\
            &\leq \|{\eta-\bar{\eta}}\|\left(1+K^2(\omega)L(\omega)\sum\limits_{s=\tau_0}^{t-1} \alpha^{t-s-1}\right)\nonumber\\
&\leq\|{\eta-\bar{\eta}}\|\left(1+\frac{K^2(\omega)L(\omega)}{1-\alpha}\right)\leq L_{G}\|{\eta-\bar{\eta}}\|.
    \end{align}

Now let us show that $H_\omega$ satisfy a global Lipschitz condition.
 \begin{eqnarray}\label{lip_H}
       \|{H_\omega(t,\xi)-H_\omega(t,\bar{\xi})}\|_{t+1,\omega}  &\le& \|{\xi-\bar{\xi}}\|+K(\omega)L(\omega)\sum\limits_{s=\tau_0}^{t-1} \alpha^{(t-s-1)}\left\{\|{\Phi_{\omega}(s,t)\xi-\Phi_{\omega}(s,t)\widetilde{\xi}}\|_{s,\omega}\right.\nonumber\\&+&\left.\|{\varphi^*_{t,\omega}(\xi)-\varphi^*_{t,\omega}(\bar{\xi})}\|_{\omega}\right\} \le \|{\xi-\bar{\xi}}\|\left(1+\frac{K^2(\omega)L(\omega)}{1-\alpha}\right)\nonumber\\&+&K(\omega)L(\omega)\sum\limits_{s=\tau_0}^{t-1} \alpha^{(t-s-1)}\|{\varphi^*_{t,\omega}(\xi)-\varphi^*_{t,\omega}(\bar{\xi})}\|_{\omega}\nonumber\\
       &\le& \|{\xi-\bar{\xi}}\|\left(1+\frac{K^2(\omega)L(\omega)}{1-\alpha}+\frac{K^3(\omega)L^2(\omega)\ell(\tau,\omega)}{(1-\alpha)(1-K(\omega)L(\omega))}\right)\nonumber\\&\leq& L_{H}\|{\xi-\bar{\xi}}\|.
    \end{eqnarray}

Based on \eqref{Lip_G} and \eqref{lip_H}, we conclude that $H_\omega(t,\cdot)$ and $G_\omega(t,\cdot)$ are Lipschitz.

which completes the proof.
\end{proof}

\subsection{Smooth conjugacy}
We  write $D_2^{j}F_\omega:\mathbb{I}\times X\to L_{j}(X)$ for the $j$-th order partial  derivative of a mapping $F_\omega:\mathbb{I}\times X\to X$ with respect to the second variable and we substitute the assumption $(H_2^0)$ by the following condition supposed to hold for $m\in\mathbb{N}$ and all $\omega\in\Omega$
\begin{enumerate}
	\item[$(H_2^m)$] there exist reals $M_j(\omega)\geq 0$ such that for any fixed $t\in \mathbb{I}$, all $x\in X$ and $0\leq j\leq m$
	\begin{align}
		\|{D_2^{j}F_{\omega}(t,x)}\|_{t+1,\omega}&\leq M_j(\omega).
	\end{align}
\end{enumerate}

Then the smooth linearization result follows.

\begin{theorem}
   \label{thmdiff}
	Assume the hypothesis $(H_1)$ and $(H_2^m)$ hold on $\mathbb{I}$ and \begin{equation}\label{est_thm2}
	M_1(\omega)K(\omega)<1-\alpha.
	\end{equation}

Then all statements of the Theorem \ref{thm_lin} hold with  $H_{\omega}{(t,\cdot)}$ being  a $C^m-$ diffeomorphism for any fixed $t\in\mathbb{I}$.
\end{theorem}

\begin{proof}

First of all we want to show that
the general solution  $\varphi_{\omega}(s,t,\cdot):X\to X$ of \eqref{deq} is $m$ times differentiable with bounded partial derivatives
$D_3^j\varphi_{\omega}:I\times I\times X\to L_{j}(X)$
for any $\eta\in X,$ fixed $s,t\in \mathbb{I}$ and $1\leq j\leq m.$

This fact is a direct consequence from the definition of the general solution. Indeed, if we abbreviate $f_{\omega}(t,x) := A_{\omega}(t)x+F_{\omega}(t,x),$ we see that $f_{\omega}(t,\cdot)$ is of class $C^m.$ Then the general solution $\varphi_{\omega}(s,t,\cdot)$ of \eqref{deq} for any $s\leq t$ can be represented as $$\varphi_{\omega}(s,t,\cdot) = f_{\omega}(t-1,\cdot)\circ...\circ f_{\omega}(s,\cdot)$$ and, by applying the chain rule, the statement follows.

 Our goal now is to show the bounded properties of the partial derivatives.

From the fact that $\varphi_{\omega}$ is the general solution of \eqref{deq} it follows that for any $\eta\in X,$  $D_3\varphi_{\omega}(s,t,\cdot)$ satisfies the differential equation
\begin{equation}\label{dvarphi}
\begin{cases}
x_{t+1} = A_{\omega}(t)x_t+D_2F_{\omega}(t,\varphi_{\omega}(s,t,\eta))x_t,\\
 x_s = I_X.\\
\end{cases}
 \end{equation}

Hence by variation of constants formula we have
\begin{equation*}
  D_3\varphi_{\omega}(s,t,\eta) = \Phi_{\omega}(t,s)+\sum\limits_{\tau = s}^{t-1}\Phi_{\omega}(t,\tau+1) D_2F_{\omega}(\tau,\varphi_{\omega}(\tau,t,\eta))D_3\varphi_{\omega}(\tau,t,\eta)
\end{equation*}
and we estimate
\begin{eqnarray*}
    \|{D_3\varphi_{\omega}(s,t,\eta)}\|_{t+1,\omega}&\leq&\sum\limits_{\tau = s}^{t-1}\|{\Phi_{\omega}(t,\tau+1)}\|_{\tau,t,\omega}\|{D_2F_{\omega}(\tau,\varphi_{\omega}(\tau,t,\eta))}\|_{\tau+1,\omega}\|{D_3\varphi_{\omega}(\tau,t,\eta)}\|_{\tau+1,\omega}\\
    &+&\|{\Phi_{\omega}(t,s)}\|_{s,t,\omega}\leq K(\omega){\alpha}^{t-s}\\&+&K(\omega)M_1(\omega)\sum\limits_{\tau = s}^{t-1}{\alpha}^{t-\tau-1}\|D_3\varphi_{\omega}(\tau,t,\eta)\|_{\tau+1,\omega}.
\end{eqnarray*}
Multiplying both sides by
 $\alpha^{s-t}$ and using Gronwall's
inequality, gives us
\begin{equation*}\displaystyle
   \alpha^{s-t} \|{D_3\varphi_{\omega}(s,t,\eta)}\|_{t+1,\omega}\leq K(\omega)\left(1+\frac{K(\omega)M_1(\omega)}{\alpha}\right)^{t-s}.
\end{equation*}
From the estimate above, we obtain
\begin{equation*}
\|D_3\varphi_{\omega}(s,t,\eta)\|_{t+1,\omega}\leq K(\omega)(\alpha+K(\omega)M_1(\omega))^{t-s} ,
\end{equation*}
and the boundedness of the first-order partial derivatives  follows directly. Furthermore, the boundedness  of higher-order partial derivatives can be established using similar arguments by induction.

Now we want to show the map $\eta\mapsto G_{\omega}{(t,\eta)}$ is  $C^m.$

This follows from the fact that $\varphi_{\omega}(s,t,\cdot):X\to X$ is $C^m,$ hypothesis $(H_2^m)$ and  from chain rule, applied to the function $\Phi_{\omega}(t,s)F_{\omega}(s,\varphi_{\omega}(s,t,\eta)).$

This allows us to conclude the desired statement and to express $D_2G_{\omega}(t,\eta)$ in explicit form, using the chain rule:
\begin{equation}\label{dH}D_2G_{\omega}(t,\eta) = I_X-\sum\limits_{s=\tau_0}^{t-1}\Phi_{\omega}(t,s+1)D_2F_{\omega}(s,\varphi_{\omega}(s,t,\eta))D_3\varphi_{\omega}(s,t,\eta).\end{equation}

\begin{eqnarray*}\label{norm_derivetiveG}
\left\|\sum\limits_{s=\tau_0}^{t-1}\Phi_{\omega}(t,s+1)D_2F_{\omega}(s,\varphi_{\omega}(s,t,\eta))D_3\varphi_{\omega}(s,t,\eta)\right\|_{t,\omega} &\leq& \frac{K(\omega)M_1(\omega)}{1-\alpha}.\end{eqnarray*}

This estimate, together with \eqref{est_thm2}, shows that the norm above is less than $1.$ Therefore, the derivative of $ G_{\omega}(t, \eta)$  is invertible everywhere.  Hence, by the Theorem~\ref{ift}, the map $ \xi \mapsto G_{\omega}(t, \xi) $ is a local  $ C^m $ diffeomorphism for any fixed $t\in\mathbb{I}$.




As the final step, since we have shown that
$G_{\omega}$
  is a homeomorphism, therefore it preserves the compactness of sets. Thus, we can apply Theorem~\ref{cacc} to
$G_{\omega},$
  which completes the proof.

\end{proof}

\section{Discrete random dynamical systems}

The aim of this section is to use the properties of Lyapunov exponents and appropriate norm construction to apply the smooth results obtained in the preceding section to discrete random dynamical systems. From now we limit to finite-dimensional Banach spaces $X$ to replicate the constructions from \cite{Arnold, Wanner_95}.

Before presenting the main results of this section, let's first introduce the necessary notation and definitions.

\begin{definition}
Let $(\Omega, \mathcal{F},\mathbb{P})$ be a probability space. We call $(\Omega, \mathcal{F},\mathbb{P},(\theta^t)_{t\in \mathbb{T}})$ a \emph{measurable dynamical system} (MDS) if the mapping $\theta^t:\Omega\to\Omega$ satisfy the following conditions
\begin{itemize}
    \item[(i)] the mapping $(t,\omega)\mapsto\theta^t\omega$ is measurable;
    \item[(ii)]the family $(\theta^t)_{t\in \mathbb{T}}$ forms a group, i.e. we have $\theta^0 = Id_{\Omega}$ and $\theta^{t+s} = \theta^t\circ\theta^s,$ for arbitrary $s,t\in \mathbb{T};$
    \item[(iii)] the mapping $\theta^t$ is $\mathbb{P}-$preserving, i.e. for arbitrary $t\in \mathbb{T}$ and $F\in \mathcal{F}$ the identity $\mathbb{P}(\theta^t(F)^{-1}) = \mathbb{P}(F)$ holds.
\end{itemize}
\end{definition}

A MDS $(\Omega, \mathcal{F},\mathbb{P},(\theta_t)_{t\in \mathbb{T}})$ is called \emph{ergodic} if every $\theta^t-$ invariant set has probability $0$ or $1,$ i.e. if for all $F\in \mathcal{F}$ satisfying $\theta^t(F)^{-1} = F$ for every $t\in \mathbb{T}$ we have either $\mathbb{P}(F) = 0$ or $\mathbb{P}(F) = 1.$
\begin{definition}

A \emph{measurable random dynamical system} (RDS)  on $X$
over some given MDS $(\Omega, \mathcal{F},\mathbb{P},(\theta^t)_{t\in \mathbb{T}})$ is a measurable mapping
$\psi : \mathbb{T} \times \Omega \times X \to X,$ $(s, \omega, x) \mapsto \psi(s, \omega, x)$ forming a \emph{cocycle} over $\theta^t,$ i.e.
 the mapping $\psi(t, \omega) := \psi(t, \omega, \cdot)$ satisfy

$$\psi(0, \omega) = I_X,$$
$$\psi(s + t, \omega) = \psi(t, \theta^s \omega) \circ \psi(s, \omega), $$
for arbitrary $s,t \in \mathbb{T},$ and $\omega\in \Omega.$
\end{definition}

The random dynamical system $\psi$ is called \emph{linear} if the mapping $\psi(t,\omega)$ is linear.

\begin{definition} Let $\Phi:\mathbb{T}\times\Omega\to X$ be a linear cocycle over an ergodic MDS $\theta^t$.
    We define the \emph{Lyapunov exponents} of $\Phi$  for all $x\in X, x\neq 0$ in the following way (provided it exists)
    $$\lambda(\omega,x) = \lim_{t\to\infty}\frac{1}{t}\|{\Phi(t,\omega)x}\|.$$
We write $-\infty<\lambda_{k(\omega)}(\omega)<...<\lambda_1(\omega)$ for the different values of $\lambda(\omega,x)$ can take on for $x\neq 0$ and call $\lambda_1(\omega)$ the \emph{top Lyapunov exponent.}

In the setting of the Multiplicative Ergodic Theorem~\ref{MET}, the sets $U_{\lambda}:=\{x:\lambda(\omega,x)\leq\lambda(\omega)\}$ are linear subspaces of $X,$ $U_i:=U_{\lambda_i}, $ forms a filtration (flag of subspaces)
$$\{0\}\subset U_{k(\omega)}\subset...\subset U_1 = X,$$
and
$$\lambda(\omega,x) = \lambda_i(\omega)\iff x\in U_i\setminus U_{i+1},\quad i=1,...,k(\omega).$$
We say that $d_i(\omega):=\dim U_i-\dim U_{i+1}$ is the \emph{multiplicity} of $\lambda_i(\omega).$
The set $$S(\theta,\Phi):=\{\lambda_i(\cdot),d_i(\cdot):i=1,...,k(\cdot)\}$$
is called the \emph{Lyapunov spectrum} of $\Phi.$
\end{definition}

From now we will focus on the \emph{discrete time random dynamical system} case with $\mathbb{T} = \mathbb{I}.$
 To achieve this, we introduce the time-one mapping
\[
\psi^{'}(\omega) := \psi(1, \omega): X \to X.
\]
By iteratively applying the cocycle property, we derive:
\[
\psi(t, \omega) =
\begin{cases}
    \psi^{'}(\theta^{t-1} \omega) \circ \cdots \circ \psi^{'}(\omega), & \text{if } t \geq 1, \\
    I_X, & \text{if } t = 0.
\end{cases}
\]

We say that the random dynamical system \( \psi \) is generated by $\psi^{'}$. It is of class \( C^k \) if and only if the mapping \( x \mapsto \psi^{'}(\omega)x \) is of class \( C^k \) for all \( \omega \in \Omega \).

To apply the results from the previous section to the random dynamical system, we first need to address the question of how the difference equation \eqref{deq} relates to the corresponding random dynamical system. The answer is provided by the following Lemma.

\begin{lemma}
    Assume there is a given arbitrary MDS $(\Omega,\mathcal{F},\mathbb{P},(\theta^t)_{t\in \mathbb{I}}),$ as well as a difference equation \eqref{deq}. Then the following statements are equivalent
    \begin{itemize}
        \item[(i)] The mapping $\psi_{\omega}$ defined as $\psi_{\omega}(t,x):=\varphi_{\omega}(t, 0,x)$ is a measurable RDS over $(\Omega,\mathcal{F},\mathbb{P},(\theta^t)_{t\in \mathbb{I}}).$
        \item[(ii)] The general solution of \eqref{deq} satisfies \begin{equation}\label{lem32}
            \varphi_{\omega}(t,\tau, x) = \varphi_{\theta^{\tau} \omega}(t-\tau, 0, x)
        \end{equation} for arbitrary $\tau,t\in \mathbb{I}$ and $x\in X.$
        \item[(iii)] The mapping $A_{\omega}(t)x+F_{\omega}(t,x)$ satisfy the identity
        \begin{equation}
            A_{\omega}(t)x+F_{\omega}(t,x) = A_{\theta^t\omega}(0)x+F_{\theta^t\omega}(0,x)
        \end{equation}
        for arbitrary $t\in \mathbb{I},$ $\omega\in\Omega$ and $x\in X.$
    \end{itemize}
\end{lemma}

\begin{proof} Let us fix $t_0, t_1, t_2\in \mathbb{I}.$
$\\(i)\Rightarrow(ii)$

From the cocycle property of RDS we know that \begin{equation}\label{a}
 \varphi_{\omega}(t_1+t_2,0,x)=\varphi_{\theta^{t_1} \omega}(t_2,0,\varphi_{\omega}(t_1,0,x)).
\end{equation}

Moreover, it is well known that the general solution of \eqref{deq} satisfies
\begin{equation}\label{b}
\varphi_{\omega}(t_2,t_0,x)=\varphi_{\omega}(t_2,t_1,\varphi_{\omega}(t_1,t_0,x)).
\end{equation}

Then
\begin{equation*}\varphi_{\omega}(t_2,t_0, x) = \varphi_{\omega}((t_2-t_0)+t_0,t_0, x)\stackrel{\eqref{a}}{=}\varphi_{\theta^{t_0} \omega}(t_2-t_0,t_0,\varphi_{\omega}(t_0,0,x)) \stackrel{\eqref{b}}{=}\varphi_{\theta^{t_0} \omega}(t_2-t_0,0,x).\end{equation*}

$(ii)\Rightarrow(i)$

In order to prove this implication it is enough to check the cocycle property.

Namely,
\begin{equation*}\varphi_{\omega}(t_1+t_2,0, x) \stackrel{\eqref{b}}{=} \varphi_{\omega}(t_1+t_2,s,\varphi_{\omega}(t_1,0,x)) \stackrel{\eqref{lem32}}{=}\varphi_{\theta^{t_1} \omega}(t_2,0,\varphi_{\omega}(t_1,0,x))\end{equation*}
and the equivalence of $(i)$ and $(ii)$ follows.

$(ii)\Rightarrow(iii)$

Let us abbreviate $f_{\omega}(t,x):=A_{\omega}(t)x+F_{\omega}(t,x)$ we see that
\begin{equation*}f_{\omega}(t, x) = f_{\omega}(t,\varphi_{\omega}(t,t,x)) = \varphi_{\omega}(t+1,t,x)\end{equation*}
and then
\begin{equation*}f_{\theta^t\omega}(0, x)  = \varphi_{\theta^t\omega} (1,0,x)\stackrel{\eqref{lem32}}{=} \varphi_{\omega}(t+1,t,x) = f_{\omega}(t,x).\end{equation*}

$(iii)\Rightarrow(ii)$

From \eqref{def_gen}, the recursive definition of the general solution, it follows that
\begin{eqnarray*}\varphi_{\omega}(t+1,s,x) &= & f_{\omega}(t,\varphi_{\omega}(t,s,x)) = f_{\theta^{t-s}\theta^s\omega}(0,\varphi_{\omega}(t,s,x))=f_{\theta^s\omega}(t-s,\varphi_{\omega}(t,s,x))\\
&=&f_{\theta^s\omega}(t-s,\varphi_{\omega}(t-s+s,s,x)) = f_{\theta^s\omega}(t-s,\varphi_{\theta^s\omega}(t-s,s,\varphi_{\omega}(s,o,x)))\\
&=& f_{\theta^s\omega}(t-s,\varphi_{\theta^s\omega}(t-s,0,x)) = \varphi_{\theta^s\omega}(t-s+1,0,x)
\end{eqnarray*}
and the lemma follows.
\end{proof}

Another necessary preliminary step is to define the new norm depending on both time and randomness.

In the settings of the Multiplicative Ergodic Theorem  \ref{MET} let $\theta^t$ is ergodic MDS, $\Phi$ be a linear RDS, $\lambda_1,...,\lambda_{k}$ be Lyapunov exponents.
Let $a>0$ denote an arbitrary, fixed, real constant such that the intervals $[\lambda_i -a, \lambda_i + a],$ $i=1,...,k,$ are disjoint. We define for arbitrary $\omega\in\Tilde{\Omega}$ and $x=x^1+...+x^{k}\in U_1(\omega)\oplus...\oplus U_{k}(\omega)$
 a new norm $|x|_{\omega} = \sqrt{|x^1|^2_{\omega}+...+|x^{k}|^2_{\omega}},$ where for $u\in U_i(\omega)$ we set $\displaystyle|u|_{\omega} = \left(\sum\limits_{t=0}^{\infty}\|\Phi(t,\omega)u\|^2 e^{-2(\lambda_i t+at)}\right)^\frac{1}{2}.$  We also assume that $|x|_{\omega} = \|x\| \text{ for all } \omega\notin \Tilde{\Omega}.$  Then (according to e.g. \cite{Arnold}, \cite{Wanner_95}) the following holds
 \begin{itemize}
     \item[(i)] $|\cdot|_{\omega}$ is a random norm on $X.$
     \item[(ii)] For every $\varepsilon >0$ there is a measurable mapping $B_{\varepsilon}:\Omega \to [1,\infty)$ such that for every $x\in X$ and almost every $\omega \in \Tilde{\Omega}$ and $t \in I$ hold
     \begin{align*}
		\frac{1}{B_{\varepsilon}(\omega)}\|x\|&\leq|x|_{\omega}\leq B_{\varepsilon}(\omega)\|x\|,&\\
		B_{\varepsilon}(\omega)e^{-\varepsilon t}&\leq B_{\varepsilon}(\theta^t\omega)\leq B_{\varepsilon}(\omega)e^{\varepsilon t}.
	\end{align*}
     \item[(iii)] For almost every $\omega \in \Tilde{\Omega},$ $t \in I$ and every $i=1,...,k(\omega), x\in U_i(\omega)$ we have
     \begin{equation}\label{random_est}
         e^{(\lambda_i-a)t} \leq \left|\Phi(t,\omega)|_{U_i(\omega)}\right|_{\omega,\theta^t\omega}\leq e^{(\lambda_i+a)t},
     \end{equation}
     where $\left|\Phi(t,\omega)|_{U_i(\omega)}\right|_{\omega,\theta^t\omega}:=\sup\{\left|\Phi(t,\omega)x\right|_{\theta^t\omega}:x\in U_i(\omega),|x|_{\omega}\leq 1\}.$
 \end{itemize}

  We consider the case, when all of Lyapunov exponents of the system are negative, i.e. $\lambda_{k}<...<\lambda_1<0$ and the constant
$a$ is chosen in such a way that
$\lambda_i+a$
remains a negative expression for all
$i=1,...,k.$

\subsection{Smooth linearization of random dynamical systems}
After establishing the connection between random dynamical systems and  difference equations, we can begin to apply the results of the previous chapter to RDS.  To do this, let $\Phi$ be a linear random dynamical system on $X$ over the ergodic MDS $(\Omega,\mathcal{F}, \mathbb{P}, (\theta^t)_{t\in \mathbb{I}}),$ satisfying conditions of  the MET and generated by the following linear equation
\begin{equation}\label{rds2}
    x_{t+1} = A(\theta^t\omega)x_t.
\end{equation}
Together with $\Phi$ we will consider
a random dynamical system $\psi$ on $X$ over the ergodic MDS $(\Omega,\mathcal{F}, \mathbb{P}, (\theta^t)_{t\in \mathbb{I}})$ generated by
\begin{equation}\label{rds1}
    x_{t+1} = A(\theta^t\omega)x_t+ F(\theta^t\omega,x_t).
\end{equation}
Assume that $\psi$ has $0$ as a fixed point, i.e. $\psi(t,\omega,0)= 0$ for all $t\in \mathbb{I},$ $\omega\in\Omega.$
Let the "nonlinear" part of $\psi$ is given by $\Psi(t,\omega,x) = \psi(t,\omega,x)-\Phi(t,\omega)x.$

First, we aim to demonstrate the validity of $(H1)$. Let $\Phi$ be a linear random dynamical system  that satisfies the Multiplicative Ergodic Theorem. Additionally, let $\Phi_{\omega}(t,s)$ represent the evolution operator of the linear random differential equation $x_{t+1} = A_{\omega}(t)x_t$. As established in \cite{Arnold} and \cite{Wanner_95}, there is a correspondence between these two objects, given by the relation $\Phi_{\omega}(t,s) = \Phi(t-s, \theta^s \omega)$.

Now, we show that in this case, hypothesis $(H1)$ holds with the norm $\|x\|_{t,\omega}:=|x|_{\theta^t\omega}$, where $K(\omega) = 1$ and $\alpha = e^{\lambda_1 + a}$.

    Indeed, from \eqref{random_est} we have
    \begin{eqnarray*}
     \left|\Phi(t,\omega)x\right|_{\theta^{t}\omega}& =& \sqrt{\left|\Phi(t,\omega)x^1\right|_{\theta^{t}\omega}+...+\left|\Phi(t,\omega)x^k\right|_{\theta^{t}\omega}} \nonumber \\
     &\leq& \sqrt{(e^{(\lambda_1+a)t}|x^1|_{\omega})^2+...+(e^{(\lambda_k+a)t}|x^k|_{\omega})^2}\nonumber\\
&\leq& e^{(\lambda_1+a)t}\sqrt{|x^1|_{\omega}^2+...+|x^k|_{\omega}^2}\nonumber\\
&\leq& e^{(\lambda_1+a)t}|x|_{\omega},
    \end{eqnarray*}
and then for any fixed $t\in \mathbb{I}$
\begin{equation}\label{analog_hyp11}
    |\Phi(t,\omega)|_{\omega,\theta^t\omega} \leq e^{(\lambda_1+a)t}.
\end{equation}
   Thus,
   \begin{equation*}
  \|{\Phi_{\omega}(t,s)}\|_{s,t,\omega} = |\Phi(t-s,\theta^s\omega)|_{\theta^s\omega,\theta^t\omega} = |\Phi(t-s,\theta^s\omega)|_{\theta^s\omega,\theta^{t-s}\theta^s\omega}
   \end{equation*}
   which, together with the~\eqref{analog_hyp11} gives us
   \begin{equation}\label{Psi_theta_S}
      \|{\Phi_{\omega}(t,s)}\|_{s,t,\omega}\leq e^{(\lambda_1+a)(t-s)}.
   \end{equation}

To satisfy the remaining hypotheses, it is necessary to impose additional conditions on the random dynamical system.

\begin{enumerate}
\item[$(H_3^0)$] Assume that  for  any $t\in \mathbb{I},$ all  $\omega\in\Omega,$ and all $x,\bar x\in X$ there exist reals $L,M\geq 0$ such that
	\begin{align}
|\Psi(t,\omega,x)|_{\theta^t\omega}&\leq M
		\label{hyp30}\\
		|\Psi(t,\omega,x)-\Psi(t,\omega,\bar x)|_{\theta^t\omega}&\leq L|x-\bar x|_{\omega}.
		\label{hyp31}
	\end{align}
\end{enumerate}

Then the following theorem establish the connection between random dynamical system $\psi$ and its linearization.

\begin{theorem}[Topological linearization of RDS]\label{thm_top}
    Assume that the linear random dynamical system $\Phi,$ generated by \eqref{rds2},  satisfies the assumptions of the MET \ref{MET} and has negative Lyapunov exponents. Assume that
   the nonlinear part $\Psi$  of a RDS $\psi,$ generated by \eqref{rds1},  satisfies $(H_3^0).$ If $L\leq \alpha,$
    then the RDS $\psi $ and $\Phi$ are topologically equivalent, i.e. there is a measurable mapping $h:\Omega \times X\to X$ with the following properties
   \begin{enumerate}
        \item[(i)]  $h(\omega) = h(\omega,\cdot)$ is a homeomorphism on $X$ with $h(\omega, 0) = 0$  for almost all $\omega\in\Omega;$
        \item[(ii)] homeomorphism maps $\omega-$orbits of $\Phi$ onto $\omega-$orbits of $\psi$ in the sense that for every $\omega\in\Omega,$ arbitrary  $t\in \mathbb{I}$ and $x\in X$
        $$h(\theta^t\omega)^{-1}\Phi(t,\omega)h(\omega)x = \psi(t,\omega, x);$$
        \item[(iii)] h has property of beeing near identity, i.e. for all $\xi,\eta\in X$
        $$|h(\omega,\xi)-\xi|_{\theta_t\omega}\leq \frac{M}{1-e^{\lambda_1+a}} \quad\text{and}\quad |h(\omega, \eta)^{-1}-\eta|_{\theta_t\omega}\leq \frac{M}{1-e^{\lambda_1+a}}.$$
    \end{enumerate}
\end{theorem}

\begin{proof}
By applying Theorem \ref{thm_lin} pointwise (for each $\omega$) to the system \eqref{rds1}, where  $\Omega$ is replaced by $\Tilde{\Omega},$ we derive the following mapping $h:\Omega\times \mathbb{I}\times X\to X:$
\begin{equation}\label{h(w)}
    \begin{cases}h(\omega,\xi) = \xi+\sum\limits_{s=0}^{t-1}\Phi(t-s,\theta^s\omega)\Psi(s,\omega,\varphi_{t,\omega}^*(\xi)(s)+\Phi(s-t,\theta^t\omega)\xi),\quad \omega\in\Tilde{\Omega},\\
   h(\omega,\xi) = 0,\quad\omega\notin\Tilde{\Omega},
   \end{cases}
\end{equation}
with the inverse given by:
\begin{equation}\label{h-1(w)}
    \begin{cases}h^{-1}(\omega,\eta) = \eta-\sum\limits_{s=0}^{t-1}\Phi(t-s,\theta^s\omega)\Psi(s,\omega,\psi(s,\omega,\eta)),\quad \omega\in\Tilde{\Omega},\\
   h(\omega,\eta) = 0,\quad\omega\notin\Tilde{\Omega},    \end{cases}
\end{equation}
ensuring that both statements
$(i)$ -
$(iii)$ hold.

Indeed, it is easy to verify that

\begin{eqnarray*}
   |h^{-1}(\omega,\eta)- \eta|_{\theta_t\omega}&\stackrel{\eqref{h-1(w)}}{\leq}&\sum\limits_{s=0}^{t-1}|\Phi(t-s,\theta^{s}\omega)|_{\theta^s\omega,\theta^t\omega}|\Psi(s,\omega,\psi(s,\omega,x))|_{\theta^t\omega}\nonumber\\
   &\stackrel{\eqref{Psi_theta_S}}{\leq}&\sum\limits_{s=0}^{t-1}e^{(\lambda_1+a)(t-s)}|\Psi(s,\omega,\psi(s,\omega,x))|_{\theta^t\omega}\stackrel{\eqref{bound_psi}}{\leq} M \sum\limits_{s=0}^{t-1}e^{(\lambda_1+a)(t-s)}\nonumber\\
   &\leq&\frac{M}{1-e^{\lambda_1+a}}.
\end{eqnarray*}
Using the same arguments we obtain:
\begin{equation*}
    |h(\omega,\xi)-\xi |_{\theta_t\omega}\leq \frac{M}{1-e^{\lambda_1+a}}.
\end{equation*}

\end{proof}

\begin{theorem}[Smooth linearization of RDS]\label{thm_smooth}
    Assume that the linear random dynamical system $\Phi$ generated by \eqref{rds2}  satisfies the assumptions of the MET \ref{MET} and has negative Lyapunov exponents. Assume
       there exist reals $M_j\geq 0$ such that
	\begin{align}\label{bound_psi}
		|D_2^{j}\Psi(\omega,x)|_{\theta^t\omega}&\leq M_j
	\end{align}
 and  \begin{equation}
	M_1(\omega)K(\omega)<1-\alpha
	\end{equation}
    for any  $t\in \mathbb{I}$, all $\omega\in\Omega,$   $x\in X$ and $0\leq j\leq m,$
 then all statements of the previous theorem hold with $h(\omega,\cdot)$ being a $C^m-$ diffeomorphism.
\end{theorem}
\begin{proof}

In order to establish that the homeomorphism \eqref{h(w)} and its inverse \eqref{h-1(w)} are of class \( C^m \), it suffices to show that both \( \Phi \) and \( \Psi \) belong to the class \( C^m \). From, e.g., \cite{Arnold} and \cite{Chu_02}, it is known that a random dynamical system is of class \( C^m \) if and only if its generator is of class \( C^m \). Hence, the desired result follows directly from the hypotheses on the operators \( A \) and \( F \).

\end{proof}

\subsection{Local smooth results}
The previous theorem's assumption that a small Lipschitz constant is necessary in practice is frequently unfeasible. Such limitations do not apply to nonlinearities that arise in real-world problems. To
address this, we focus on a local result using a cut-off technique. This method allows us to
limit the effect of the nonlinearity by "cutting off" its influence when it becomes too large.
The strategy is supported by initial results from the next lemmas, which offer important prerequisites for the efficacy of this truncation technique. These lemmas allow us to handle the issue locally without requiring a Lipschitz constant that is small globally.

\begin{lemma}\label{wan_1}\cite[Lemma 4.1, p. 257]{Wanner_95}
   Let $(\hat{\Omega},\hat{\cF})$ denote an arbitrary measurable space and  $f:\hat{\Omega}\times X\to X$ a measurable mapping such that $f(\hat{\Omega},\cdot)$ is continuous for all $\omega\in\hat{\Omega}$ with $f(\omega,0) = 0$  and
   \begin{equation}\label{lem_cut_of1_1}
     \lim_{(x,y)\to(0,0)}\frac{\norm{f(\omega,x)-f(\omega,y)}}{\norm{x-y}} = 0.
   \end{equation}
    Let $\Tilde{L}>0$ be an arbitrary constant. Then there are measurable mappings $\sigma:\hat{\Omega}\to \mathbb{R}^+$ and $\Tilde{f}:\hat{\Omega}\times X\to X$  such that the following holds:
   \begin{itemize}
    \item[(i)] if we define the random neighborhood of $0$ by
    \begin{equation}\label{U}
    U(\omega):=\{x\in X:\|x\|< \sigma(\omega)\},\end{equation}
    then the identity $f(\omega,x) = \Tilde{f}(\omega,x)$ holds for all $\omega\in\hat{\Omega}$ and $x\in U(\omega);$
    \item[(ii)] for arbitrary $\omega\in\hat{\Omega}$ and $x,y\in X$ we have
    $$|\Tilde{f}(\omega,x)-\Tilde{f}(\omega,y)|_{\theta^1\omega}\leq L(\omega)|x-y|_{\omega}\quad\text{ and }\quad |\Tilde{f}(\omega,x)|_{\theta^1\omega}\leq 1.$$

   \end{itemize}
\end{lemma}

\begin{lemma}\cite[Prop. 4.1, p. 258]{Wanner_95}
    Assume we are given a discrete-time random dynamical system $\psi$ on $X$ over the ergodic metric dynamical system $(\Omega,\cF, \mathbb{P}, (\theta^t)_{t\in \mathbb{I}}),$ with fixed point $0.$ Assume further that there is a linear random dynamical system $\Phi,$ and that the nonlinear part $\Psi$ defined via $\Psi(t,\omega,x) = \psi(t,\omega,x) - \Phi(t,\omega)x$ satisfies
    $$\lim_{(x,y)\to(0,0)}\frac{\norm{\Psi(1,\omega,x)-\Psi(1,\omega,y)}}{\norm{x-y}} = 0 $$
for arbitrary $\omega\in \hat{\Omega}.$
Then there is a  measurable mapping $\sigma:\hat{\Omega}\to\mathbb{R}^+$ and a discrete-time random dynamical system $\Tilde{\psi}(t,\omega,x) = \Phi(t,\omega)x+\Tilde{\Psi}(t,\omega,x)$ over $(\Omega,\cF, \mathbb{P}, (\theta_t)_{t\in \mathbb{I}})$ satisfying $(H_3^0),$ as well as
\begin{equation}\label{loc_lin_eq}
  \psi(t,\omega,x) = \Tilde{\psi}(t,\omega,x)
\end{equation}
for $\omega\in\hat{\Omega},$ $x\in U(\omega)$ and $t\in \mathbb{I}_{\max}(\omega,x):=\{0,t_{\max}(\omega,x)\},$ where $U(\omega)$ is defined in \eqref{U} and
\begin{equation}\label{t_max}
t_{\max}(\omega,x):=\sup\{t\in \mathbb{I}: \psi(\tau,\omega,x)\in U(\theta_{\tau}\omega) \fall 0\leq\tau\leq t\}\geq 0. \end{equation}

Moreover, for every $\omega\in\hat{\Omega}$ we have $\displaystyle\lim_{x\to 0} t_{\max}(\omega,x) = \infty.$ In particular, this implies that for every choice of $\omega\in\hat{\Omega}$ and $\tau_0\in \mathbb{I}$ the existence of a neighborhood $U_{\tau_0}(\omega)\subset X$ of the origin such that the equality \eqref{loc_lin_eq} holds for all $t\in \{1,...,\tau_o\}$ and $x\in U_{\tau_0}(\omega).$

\end{lemma}

Now, we are able to formulate a local $C^m$  result  as well.

\begin{theorem}
  [Local topological linearization of RDS]\label{thm_local}
Let $\psi$ be a  given  random dynamical system on $X$ over the ergodic metric dynamical system $(\Omega,\cF, \mathbb{P}, (\theta^t)_{t\in \mathbb{I}})$ with fixed point $0,$  generated by the random differential equation \eqref{rds1}. And and $\Phi$ be a linear random dynamical system  on $X,$ generated by the random differential equation \eqref{rds2}.
If $\Phi$ satisfies MET, $F(\omega,0) = 0$ and $$\lim_{(x,y)\to(0,0)}\frac{\norm{F(\omega,x)-F(\omega,y)}}{\norm{x-y}} = 0 $$
for all $\omega\in \Omega,$
then the  newly defined  local random dynamical system $\Tilde{\psi}$ and $\Phi$ are topologically equivalent on $\mathbb{I}_{\max}(\omega,x):=(0,t_{\max})$ in the sense of Theorem \ref{thm_top} for all $x\in U(\omega)$ and $\omega\in\hat{\Omega},$ where $t_{\max}$ defined in \eqref{t_max} and $U(\omega)$ defined in \eqref{U}.
\end{theorem}
\begin{theorem}[Local smooth linearization of RDS]\label{thm_local_smooth}
    If the assumptions of Theorem~\eqref{thm_local} hold and additionally there exist reals $M_j\geq 0$ such that all partial derivatives $D_2^{j}F:\hat{\Omega}\tm U(\omega)\to L_j(X)$,  satisfy
	\begin{align}
		|D_2^{j}F(\omega,x)|_{\theta_t\omega}&\leq M_j\quad\text{ for any  $t\in \mathbb{I^{'}}_{max}$, all $\omega\in\Omega,$   $x\in U(\omega)$ and }0\leq j\leq m,
	\end{align}
 then random dynamical system $\Tilde{\psi}$ and $\Phi$ are $C^m-$topologically equivalent  on $\mathbb{I}_{\max}.$
\end{theorem}
\begin{remark}
    The proofs of both above theorems \ref{thm_local} and \ref{thm_local_smooth} follow from the application of the Theorems \ref{thm_top} and \ref{thm_smooth} respectively to the new random dynamical system $\Tilde{\psi}$ with $X$ restricted to $U(\omega)$ and $\mathbb{I}$ restricted to $\mathbb{I}_{\max}(\omega,x).$
\end{remark}

\section*{Funding}

This work is supported by the Austrian Science Fund
(FWF): DOC 78.

\bibliographystyle{abbrv}
\bibliography{arxiv}

\setcounter{section}{0}
\renewcommand{\thesection}{\Alph{section}}
\section{Appendix}
This appendix contains all the essential results and lemmas used throughout the paper to support the proofs of our results. It is organized into two logical parts: deterministic and random, allowing readers to easily reference the material relevant to each type of system discussed in the paper.

In the deterministic part, we begin with a discrete version of Gronwall's Lemma, which is fundamental for estimating the growth of solutions of difference equations.
\begin{lemma}(Gronwall's lemma)\cite[Lem. 2.1, p. 504]{Aulbach_98}
Let $\kappa\in \mathbb{Z}$ and $b\geq 0$ be given. If the mappings $a, c: \mathbb{Z}\to\mathbf{R}^{+}_{0}$ satisfy the inequality $$a(k)\leq c(k)+b\sum\limits_{i=\kappa}^{k-1} a(i)$$
for all $k\geq\kappa,$ then for all $k\geq\kappa$ we get the estimate
$$a(k)\leq(1+b)^{\kappa - 1}c(\kappa)+ \sum\limits_{i=\kappa+1}^{k} (1+b)^{k - 1}(c(i) - c(i-1)).$$
\end{lemma}

\begin{theorem}\label{var_const}(Variation of Constants Formula)\cite[ p. 12]{Wanner_diss} Consider the system:
\begin{equation}\label{1}
   x_{t+1} = A_{\omega}(t)x_t + f(t,\omega),
\end{equation}
where \( A: \mathbb{I}\times\Omega \to L(X) \) and \( f: I \times\Omega \to X \).
Then the general solution \( \phi(t,\tau,\omega,\xi) \) of \eqref{1} is given by the variation of constants formula:
\begin{equation*}
\phi(t,\tau,\omega,\xi) =
\begin{cases}
    \xi, & \text{for } t = \tau, \\
    \Phi_{\omega}(t, \tau)\xi + \sum\limits_{i = \tau}^{t - 1} \Phi_{\omega}(t, i + 1)f(i, \omega), & \text{for } t > \tau.
\end{cases}
\end{equation*}
which is valid for all \( (t,s,\xi) \in I \times I \times X \).
\end{theorem}

Now we present the well-known fixed point theorem of Banach and some related results. From now on, let $X, Y$ be Banach spaces.

\begin{theorem}\label{Fixed point}(Banach's Fixed Point Theorem)\cite[Thm. B.1, p. 114]{Wanner_96}

Suppose $\cT : X \to X$ a contraction, i. e. there exists $c\in [0,1)$ such that

$$\norm{\cT(x)-\cT(\bar{x})}\leq c\norm{x-\bar{x}}$$
for all $x, \bar{x} \in X.$
 Then $\cT$  has exactly one fixed point $\varphi \in X.$
\end{theorem}

\begin{theorem}\label{uniform contraction}(Uniform Contraction Principle) \cite[Thm. 2.2, p.25]{Hale}

Suppose $U,$ $V$ be open sets in  $X,$ $Y$ respectively, let $\bar{U}$ be the closure of $U.$ Let
 $\cT : \bar{U}\times V \to \bar{U}$ a uniform contraction on $\bar{U}$, i. e. assume that there is a $c\in [0,1)$ such that for all $x,\bar{x}\in \bar{U}$ and $y\in V$ we have

$$\norm{\cT(x,y)-\cT(\bar{x},y)}\leq c\norm{x-\bar{x}}.$$
 Let $\varphi(y) $ be the unique fixed point of $\cT(\cdot,y)$ in $\bar{U}.$ If $\cT\in C^k( \bar{U}\times V,X),$ $0\leq k<\infty,$ then $\varphi(\cdot) \in C^k(V,X).$

\end{theorem}
\begin{theorem}\label{ift}(Local Inverse Mapping Theorem) \cite[Thm. 4.F, p. 172]{zeidler}

Let  $f:U(x_0)\subseteq X\to Y$ be a $C^m$ mapping for $m\in\mathbb{N}$. Then $f$ is a local $C^m-$ diffeomorphism at $x_0$ if and only if $Df(x_0)$ is bijective.
\end{theorem}

\begin{theorem}\label{cacc}(Global Inverse Mapping Theorem)\label{hadamard} \cite[Thm. 4.G, p. 174]{zeidler}

Let $f:X\to Y$ be a  local $C^m-$ diffeomorphism for $m\in\mathbb{N}$ at every point of $X.$  Then $f$ is a $C^m-$ diffeomorphism if and only if $f$ is proper.
\end{theorem}

Now let us turn to a key result in random dynamical systems: the fundamental Oseledets's Multiplicative Ergodic Theorem (MET).

In the following we will write $\log^+(x) = \max\{0, \log(x)\}$ for any $x\in X .$

\begin{theorem}[MET, half line, discrete time]\cite[Thm. 3.4.1, p. 134]{Arnold}\label{MET}

Let $\Phi$ be a linear cocycle over an ergodic MDS $(\Omega, \mathcal{F}, \mathbb{P}, (\theta^n)_{n \in \mathbb{Z}})$.	

Let
\[
\Phi(n, \omega) =
\begin{cases}
A(\theta^{n-1}\omega) \cdots A(\omega), & n > 0, \\
I_X, & n = 0, \\
A(\theta^n\omega)^{-1} \cdots A(\theta^{-1}\omega)^{-1}, & n < 0,
\end{cases}
\]
be generated by $A : \Omega \to Gl(d, \mathbb{R})$ and assume
\[
\log^+ \|A(\cdot)\| \in L^1(\Omega, \mathcal{F}, \mathbb{P}) \quad \text{and} \quad \log^+ \|A^{-1}(\cdot)\| \in L^1(\Omega, \mathcal{F}, \mathbb{P}). \tag{1.12}
\]	

 Then there exist a forward invariant set $\Tilde{\Omega}\in\cF$ of full measure such that for every $\omega\in\Tilde{\Omega}$ the following statements hold:
 \begin{enumerate}
    \item[(i)] There exist $k(\omega)$ numbers $\lambda_1(\omega) >...> \lambda_{k(\omega)}(\omega)$ and the invariant splitting $\mathbb{R}^d=U_1(\omega)\oplus...\oplus U_{k(\omega)}(\omega)$ such that
    \begin{align*}
        k(\theta\omega) &{=} k(\omega),\\
       \lambda_i(\theta\omega) &{=} \lambda_i(\omega) \fall i\in\{1,2,...,k(\omega)\},\\
       d_i(\theta\omega) &{=} d_i(\omega) \fall i\in\{1,2,...,k(\omega)\},
    \end{align*}
where $d_i(\omega):=\dim U_i(\omega).$
    \item[(ii)] Put $V_{k(\omega)+1}(\omega):=0$ and for $i=1,..., k(\omega)$
    $$V_i(\omega):=U_i(\omega)\oplus...\oplus U_{k(\omega)}(\omega)$$ so that
    $$V_{k(\omega)}(\omega)\subset ...\subset V_{i}(\omega)\subset ...\subset V_{1}(\omega)$$
    defines a filtration on $\mathbb{R}^d.$ Then for each $x\in \mathbb{R}^d\setminus \{0\}$ the Lyapunov exponent
    $$\lambda(\omega,x) := \lim_{n\to\infty}\frac{1}{n}\log\norm{\Phi(n,\omega)x}$$
    exists and
    $$\lambda(\omega,x) = \lambda_i(\omega)\iff x\in V_i(\omega)\setminus V_{i+1}(\omega),$$
    equivalently
    $$V_i(\omega) = \{x\in X:\lambda(\omega,x) \leq \lambda_i(\omega)\}.$$
    \item[(iii)] For all $x\in \mathbb{R}^d\setminus \{0\}$
    $$\lambda(\theta\omega,A(\omega)x) = \lambda(\omega,x),$$
    whence
    $$A(\omega) V_i(\omega)\subset V_i(\theta\omega)\fall i\in\{1,2,...,k(\omega)\}.$$
 \end{enumerate}
 \end{theorem}

\end{document}